\documentclass[11pt]{article}
\usepackage{amsmath,amssymb,amsthm,subfigure,float,url}
\usepackage[applemac]{inputenc}

  \usepackage[colorlinks=true]{hyperref}

\usepackage[dvips]{graphicx}

\newcommand{\R} {\mathbb{R}}             %
  
 \newcommand{\findemo}{\hfill $\Box$}

\newtheorem{Th}{Theorem}[section]
\newtheorem{prop}{Proposition}[section]

\newtheorem{lem}{Lemma}[section]

\newtheorem{rem}{Remark}[section]

\renewcommand{\geq}{\geqslant}
\renewcommand{\leq}{\leqslant}

\title{Smoothing Methods for Nonlinear Complementarity Problems}

\author{Mounir Haddou  and Patrick Maheux
\\
F\'ed\'eration Denis Poisson, D\'epartement de Math\'ematiques MAPMO, 
\\Universit\'e d'Orl\'eans, F- 45067 Orl\'eans, France}

\begin{document}
    
\maketitle

\begin{abstract}
In this paper, we present a new smoothing approach to solve general  nonlinear complementarity problems.  Under  the $P_0$ condition on the original problems, we prove some existence and convergence results . We also present an error estimate under a new and general monotonicity condition. The numerical tests conﬁrm the efﬁciency of our proposed methods. 
\end{abstract}

\noindent


\noindent

\section{Introduction}\label{intro}
\setcounter{equation}{0}

Consider the nonlinear complementarity problem (NCP), which is to find 
a solution of the system :
\begin{equation}\label{complem}
x\geq 0, F(x)\geq 0\quad {\rm and} \quad  x^{\top} F(x)=0,
\end{equation}
where  $F:\R^n\longrightarrow \R^n$ is a continuous function that satisfies some additional assumptions to be precise later.\\

This problem has a number of important applications in operations research, economic equilibrium 
problems and in the engineering sciences \cite{FP} . It has been extensively studied and the number of proposed solution methods  is enormous eg  (\cite{FMP}  and references therein).
There are almost three  different classes of methods: equation-based methods (smoothing) , merit functions and projection-type methods.\\
Our goal in this paper is to present  new and very simple  smoothing and approximation schemes to solve NCPs and to produce   efficient  numerical methods. These  functions are based on penalty functions for convex programs.\\


Almost all the solution methods consider at least the following  important and standard condition on the mapping $F$ ( monotonicity)
: for any $x, y\geq 0$,
\begin{equation}\label{complem}
( x-y)^{\top}(F(x)-F(y))\geq 0.
\end{equation}
We will assume that 
$$F {\rm\  is\  a\quad } (P_0){\rm -function }\leqno{(H0)}$$
to prove the convergence of our approach. This assumption is weaker than monotonicity.\\
We recall the following definitions of $(P_0)$-  and $(P)$ -functions. 
We say that     $F:\R^n\rightarrow \R^n$ is a $(P_0)$-function if 
$$
\max_{i: x_i\neq y_i} (x-y)_i(F_i(x)-F_i(y))\geq 0,\quad x,y\in \R^n,
$$
and $F$ is a $(P)$-function if 
$$
\max_{i: x_i\neq y_i} (x-y)_i(F_i(x)-F_i(y))> 0 ,\quad x,y\in \R^n.
$$
We start with an easy result. We define   component-wise  the function
$F_{min}(x):=\min(x,F(x))$ with   $F_{{min},i}(x)=\min(x_i,F_i(x))$ for any $i:1...n$.

This function possesses the same properties as $F$
\begin{lem}
Assume that $F$ is a $(P_0)$ (respectively $(P)$)-function then $F_{min}(x):=\min (x,F(x))$ is also 
$(P_0)$ (respectively $(P)$)-function.
\end{lem}
\proof
Assume that  $F$ is a $(P_0)$-function. For  any $x,y\in\R^n$, there exists $i $ such that $x_i\neq y_i$ and $(x_i-y_i)(F_i(x)-F_i(y))\geq 0$.  We can assume $x_i>y_i$. Then 
$
F_i(x)\geq F_i(y)
$.
So
$
F_i(x)\geq \min(y_i,F_i(y))
$ and $x_i\geq \min(y_i,F_i(y))$.
Then 
$
 \min(x_i,F_i(x))
\geq \min(y_i,F_i(y))
$.  Thus 
$(x_i-y_i)(F_{min,i}(x)-F_{min,i}(y))\geq 0
$. 
Then $F_{min}$ is a  $(P_0)$-function.
The proof is analogue if $F$ is a $(P)$-function.
\findemo
\\

An other ssumption that will be useful in our approach is that the   solution set is compact
$${\mathcal Z}:=\{x\geq 0, F(x)\geq 0,  x^{\top} F(x)=0\} {\rm\  is\  nonempty\  and\  compact. } \leqno{(H1)}$$

\begin{rem}
Under some sufficient conditions, the assumption (H1) is satisfied.  Note that this set may be empty, for instance if $-F(x)>0$ for any $x\in \R^n$ even when $F$ is continuous and monotone. Such counter-example is easy to show,  for example :$F(x)=\frac{-1}{x+1}$ for $x\geq 0$ and $F(x)=-1$ if $x\leq 0$.
\end{rem}
 We give in the following lemma an  example of sufficient condition on the mapping $F$ to  insure $(H1)$. 
 \\
 
\begin{lem}
Assume that $F$ is continuous and monotone on $\R^n$. Moreover, we assume that:
\begin{enumerate}
\item
There exits $y\in\R^n$ with
$F(y)>0$ .
\item
There exist constants  $c, M>0$ such that
for any $x$, $\vert x\vert_1 \geq M$,
$\vert F(x)\vert \leq c \vert x\vert_1$ with $c<m(F(y))/\vert y\vert $ where 
$m(F(y))=\min_i F_i(y)$.
\end{enumerate}

Let $\varepsilon>0$.
Then ${\mathcal Z}_{\varepsilon}:=\{x\geq 0, F(x)\geq 0,  x^{\top} F(x)
\leq \varepsilon
\}$ is compact (may be empty).
\end{lem}

\proof 
Since $F$ is continuous,  the set ${\mathcal Z}_{\varepsilon}$ is closed. To show the compactness, it is enough to show the boundedness of ${\mathcal Z}$. The monotonicity property implies for any $x\in {\mathcal Z}_{\varepsilon}$:
$$
 x^{\top} F(y)\leq  x^{\top} F(x)-y^{\top} F(x)+y^{\top} F(y)
\leq
 \vert y\vert \vert F(x)\vert +\vert y\vert \vert F(y)\vert +{\varepsilon}.
$$
Let $x\in  {\mathcal Z}_{\varepsilon}$ and $\vert x\vert_1 \geq M$,
then
$$
m(F(y)) \vert x\vert_1\leq 
\sum_{i=1}^n x_i F_i(y)\leq  
c  \vert y\vert  \vert x\vert_1
+\vert y\vert \vert F(y)\vert
+{\varepsilon}.
$$
Thus 
$$
(m(F(y))- c  \vert y\vert  )
\vert x\vert_1
\leq \vert y\vert \vert F(y)\vert +{\varepsilon}.
$$

Since $\kappa := m(F(y)) -c  \vert y\vert>0$, we get
$$
\vert x\vert_1
\leq (\vert y\vert \vert F(y)\vert +{\varepsilon})\vert/\kappa.
$$
Thus 
$$
\vert x\vert_1
\leq 
\max( M,(\vert y\vert \vert F(y)\vert +{\varepsilon})/\kappa).
$$

So ${\mathcal Z}_{\varepsilon}$ is bounded, hence compact.
\findemo
\\
\begin{rem} All   continuous monotone bounded function $F$  satisfying  the  condition (1) also satisfies condition (2)  of the lemma.
Indeed there exists $R>0$ such that for any $M>0$ and any $x$ such that $\vert x\vert_1 \geq M$,
$$
\vert F(x)\vert \leq R\leq \frac{R}{M} \vert x\vert_1 
$$
Let $c:=\frac{R}{M}$, it enough to choose $M$ large enough such that $c< m(F(y))/\vert y\vert $.
\\

This condition (2) allows us to consider a family of functions  $F$ satisfying some sub-linear growth  at infinity.
\end{rem}


The organization of the paper is as follows. In Section 2, we define the smoothing functions  and the approximation technique. In Section 3, we  give a detailed discussion of the properties of the smoothing function and the approximation scheme. Section 4 is devoted to the proof of convergence and the error estimate. Numerical examples and results will be reported  in the last section. 

\section{The smoothing  functions}

We start our discussion by introducing the function $\theta$ with the following properties (See \cite{ACH,Had}).
Let $\theta:\R\rightarrow (-\infty,1)$ be an increasing continuous function such that
$\theta(t)<0$ if $t<0$, $\theta(0)=0$ and $\theta(+\infty)=1$.
For instance $\theta^{(1)}(t)=\frac{t}{t+1}, t\geq 0$ and  $\theta^{(1)}(t)=t$ if  $t<0$,
$\theta^{(2)}(t)=1-e^{-t},\, t\in \R$. 

This function "detects" if $t=0$ or $t>0$ i.e. if $t\geq 0$ in a "continuous way".
The ideal function will be the function $\theta^{(0)}(-\infty)=-\infty$, $\theta^{(0)}(t)<0$ if $t<0$, $\theta^{(0)}(0)=0$ and 
$\theta^{(0)}(t)=1$ for $t>0$. But doing so, at least a discontinuity at $t=0$ is introduced.  We smooth this ideal function $\theta^{(0)}$ by introducing  $\theta_r(t)=\theta(t/r)$ for $r>0$.
\\
So that,  $ \theta_r(0)=0$,  $\forall r>0$, $\lim_{r\rightarrow 0^+} \theta_r(t)=1$  for
all $t>0$ and $\lim_{r\rightarrow 0^+} \theta_r(t)=\inf  \theta=\theta(-\infty)<0$ if $t<0$. So   $\lim_{r\rightarrow 0^+} \theta_r$ behaves essentially as   $\theta^{(0)}$. Moreover, note that the  function  $\theta_r$corresponding differentiates quantitatively the positive values of t: if $0<t_1<t_2$ then $0<\theta_r(t_1)<\theta_r(t_2)$ and conversely.
\\

Now, let's  consider the following equation on the one-dimensional case. Let $s,t\in \R^+$ be such that
\begin{equation}\label{theta}
\theta_r(s)+\theta_r(t)=1.
\end{equation}
For instance, let's take $\theta^{(1)}$.
The equality (\ref{theta})  is  then equivalent to 
$$st=r^2$$
So, when $r$ tends to $0$, we simply get   $st=0$. This limit case applied with $s=x\in\R^+$ and $t=F(x)\in\R^+$ gives  our relation $xF(x)=0$. Our approximation is 
$x^{(r)}F(x^{(r)})=r^2$.
So, for  general $\theta$, the aim of this paper is to produce, for each $r\in(0,r_0)$, a solution $x=x^{(r)}\in \R^n$ with $x^{(r)}\geq 0$ such that 
$F(x^{(r)})\geq 0$ and  
\begin{equation}\label{thetadim}
\theta_r(x)+\theta_r(F(x))=1.
\end{equation}
and to show the compactness of the set $\{x^{(r)}, r\in(0,r_0)\}$. Hence by taking a subsequence of $x^{(r)}$, we expect to converge to a solution of $xF(x)=0$.
The equation just above  has to be interpreted, in the multidimensional case, as
$$
\theta_r(x_i^{(r)})+\theta_r(F_i(x^{(r)}))=1, \qquad i:1...n.
$$

Note that the relation (\ref{thetadim})  is symmetric in $x$ and $F(x)$ and it can be seen as a fixed point problem for the function
$F_{r,\theta}(x)$ defined just 
below.
 Indeed,
(\ref{thetadim}) is equivalent to
$$
x= \theta_r^{-1}\left(1-\theta_r(F(x))\right)=r \theta^{-1}\left(1-\theta(F(x)/r)\right)
= :
F_{r,\theta}(x)
$$
and also, by  symmetry of the equation  (\ref{theta}) (if any), we have the relations:
$$
F(x)= \theta_r^{-1}\left(1-\theta_r(x)\right)=r \theta^{-1}\left(1-\theta(x/r)\right).
$$

This fact may be of some interest for numerical methods.
The  speed of convergence to a solution can be compared  for different choices of $\theta$. 
\\

Now, we  propose another way to approximate a  solution of the (NC) problem.
\\

Let $\psi_r(t)=1-\theta_r(t)$.
The relation (\ref{thetadim})  is equivalent to
$$
\psi_r(x)+\psi_r\left(F(x)\right)=1=\psi_r(0).
$$
Hence,  the relation can be written as
$$
\psi_r^{-1}[\psi_r(x)+\psi_r(F(x))]=0.
$$
Let $\psi=\psi_1=1-\theta$. thus, we have
$$
r\psi^{-1}\left[\psi\left(\frac{x}{r}\right)+\psi\left(\frac{F(x)}{r}\right)\right]=0.
$$
For the sequel, we set 
$$
G_r(x,y):=
r\psi^{-1}\left[\psi\left(\frac{x}{r}\right)+\psi\left(\frac{y}{r}\right)\right].
$$

First, we characterize solutions $(x,y)$  of $G_r(x,y)=0$
when $\theta$ satisfies some conditions independent of $F$.
\\

Let $0<a<1$. We say that $\theta$ satisfies 
$(H_a)$ if there exists   $s_a>0$ such that, for all  $s\geq s_a$,
$$
\frac{1}{2}+\frac{1}{2}\theta(as)\leq \theta(s).
$$
This condition is equivalent to
$$
 \psi(s)\leq \frac{1}{2}  \psi(as),\quad s\geq s_a.
$$

{\bf 
Two  examples}:

\begin{enumerate}
\item
Let  $\theta^{(1)}(t)=\frac{t}{t+1}, t\geq 0$ and  $\theta^{(1)}(t)=t$ if  $t<0$ (But the case  $t<0$ is not useful in this discussion).  Then
$\psi^{(1)}(t)=\frac{1}{t+1}$ if $t\geq 0$ and $\psi^{(1)}(t)=1-t$ if $t<0$.
The condition $(H_a)$ is only satisfied with $0<a<1/2$ and $s_a\geq \frac{1}{1-2a}$.
\item
Let $\theta^{(2)}(t)=1-e^{-t},\, t\in \R$. 
Then $\psi^{(2)}(t)=e^{-t}$  satisfies the condition $(H_a)$ for any $0<a<1$ with 
$s_a=\frac{\ln 2}{1-a}$.
\end{enumerate}

Note that these functions $\psi$ do not satisfies the condition $(H_a)$ in the same range for $a$. This has some consequence for the limite of $G_r(s,t)$ as $r$ goes to zero (See Th.\ref{allha} and  Example 1 after the proof of Th.\ref{zero}).
\\


\begin{Th}\label{zero}
Assume that for some $0<a<1$, the condition  $(H_a)$ is satisfied for $\theta$. Let $s,t\in \R$.
 The two following statements are equivalent
\begin{enumerate}
\item
$\lim_{r\rightarrow 0}G_r(s,t)=0$ 
\item
$\min(s,t)=0$.
\end{enumerate}
\end{Th}

The statement $\min(s,t)=0$ is equivalent to
$s=0\leq t$ or $t=0\leq s$. This is the reason why we expect   a solution of the (NC) problem
with $s=x_i$ and $t=F_i(x)$ by considering the function $G_r$.
\\

\proof
We asserts that
$
G_r(s,t)\leq \min(s,t)
$ for any $r>0$ and $s,t\in\R$.
Indeed,  assume that $s\leq t$ i.e. $s=\min(s,t)$.
 Since $\psi \geq 0$,
$$
\psi(s/r)\leq  \psi(s/r)+\psi(t/r)
$$
So, by the fact that $\psi$ is non-increasing,
$$
s/r\geq  \psi^{-1}\left( \psi(s/r)+\psi(t/r)\right).
$$
Then
$$G_r(s,t)\leq s.$$
The assertion is proved.
We now prove (1) implies (2). Let $s,t\in \R$,
we have assumed that
$\lim_{r\rightarrow 0}G_r(s,t)=0$ thus $ \min(s,t)\geq 0$.
By symmetry, we can suppose  that 
$s= \min(s,t)$.  We deduce the result by contradiction.
Assume that $s>0$. Since $\psi$ is non-increasing,
$$
 \psi(s/r)+\psi(t/r)\leq 2 \psi(s/r).
 $$
 For $r$ small enough,
 $2 \psi(s/r)\leq  \psi(as/r)$ because $s/r$ goes to infinity.
 Whence
 $$
 \psi(s/r)+\psi(t/r)\leq  \psi(as/r).
 $$
 Again since $\psi$ is non-increasing,
 $$
as/r\leq  \psi^{-1}\left( \psi(s/r)+\psi(t/r)\right).
$$
or equivalently ($r$ small enough),
$$
s\leq a^{-1}G_r(s,t).
$$
By assumption, $\lim_{r\rightarrow 0}G_r(s,t)=0$ then $s\leq 0$.  Contradiction. So $s=0$
and the implication is proved.
\\

We now prove the converse.
Assume 
$ s=\min(s,t)= 0$. Then 
$ G_r(s,t)=r\psi^{-1}\left( 1+\psi(t/r)\right)$ due to  $\psi(0)=1$.
If $t=0$, $ \lim_{r\rightarrow 0}G_r(s,t)=\lim_{r\rightarrow 0}r\psi^{-1}(2)=0.$
If $t>0$,
$\lim_{r\rightarrow 0}\psi(t/r)=1-\lim_{r\rightarrow 0}\theta(t/r)=0
$. Thus
$ \lim_{r\rightarrow 0}G_r(s,t)=\lim_{r\rightarrow 0}r\psi^{-1}(1)=0$ by continuity of $\psi^{-1}$.
The proof is completed.
\\

{\bf Two examples}:
\begin{enumerate}
\item
Let  $\theta^{(1)}(t)=\frac{t}{t+1}, t\geq 0$ and  $\theta^{(1)}(t)=t$ if  $t<0$.  Then
$\psi^{(1)}(t)=\frac{1}{t+1}$ if $t\geq 0$ and $\psi^{(1)}(t)=1-t$ if $t<0$. For $s>0$ and $t>0$ such that $\frac{1}{s}+\frac{1}{t}\leq \frac{1}{r}$, then
$$
 G_{1,r}(s,t)=\frac{st-r^2}{s+t+2r}.
 $$
 Note that the denominator is not zero when $s,t$ are positive even in the case $s=t=0$. This is interesting fact for numerical simulation.
\\

In that case
  $\lim_{r\rightarrow 0} G_{1,r}(s,t)=\frac{st}{s+t}$.
  Note that this is not   $\min(s,t)$. We can easily prove that
    $\lim_{r\rightarrow 0} G_{1,r}(s,t)=0$ if $s=0$ and $t>0$ or $t=0$ and $s>0$.
    \\

 If $t,s> 0$, the derivative in $r$ of $ G_{1,r}(s,t)$ is
$$
\frac{-2br-2r^2-2a}{(b+2r)^2}\leq 0
$$
 with $a=st$ and $b=s+t$.
 So $ G_{1,r}(s,t)$ is non-increasing in $r$ for fixed $s,t>0$. Since  $G_{1,r}(s,t)\leq \min(s,t)$
 then   $\lim_{r\rightarrow 0} G_{1,r}(s,t)$ always exists.  

\item
Let $\theta^{(2)}(t)=1-e^{-t},\, t\in \R$. 
Then $\psi^{(2)}(t)=e^{-t}$ and 
$$
 G_{2,r}(s,t)=-r\log (e^{-s/r}+e^{-t/r}).
 $$
for any $s,t\in \R$,
$$
\lim_{r\rightarrow 0} G_{2,r}(s,t)=\min(s,t).
$$
Indeed,  if $s=\min(s,t)$ then  $-r\log 2+s\leq G_{2,r}(s,t)$ because
$ e^{-s/r}+e^{-t/r}\leq 2  e^{-s/r}$.
Thus
$$
-r\log 2+\min(s,t)\leq  G_{2,r}(s,t)\leq \min(s,t).
$$
So, we deduce the expected limit. The assertion of Th.\ref{zero} is clealy satisfied.\\
\end{enumerate}
 
It is an easy exercice to show that these  two functions $ G_{1,r} $ and  $G_{2,r}$ and there limit function $G_{1,0}, G_{2,0}$  are concave functions on $(\R^2)^+$.
(Pb:  Is  it always true for any $\theta$?)
\\

Now, we focus on the case where  $\theta$ satisfies  $(H_a)$  for all $a\in (0,1)$.
\begin{Th}\label{allha}
Assume that $\theta$ satisfies  $(H_a)$ is for all $a\in (0,1)$.
Then for any $s,t>0$,
$$
\lim_{r\rightarrow 0}G_r(s,t)=\min(s,t).
$$
\end{Th}

\proof
We use essentially the arguments of Th.\ref{zero}.
We have seen in the course of the proof  that
$
G_r(s,t)\leq min(s,t)
$ for any $r>0$, $s,t\in \R$ and any $\psi$.  
\\

Now  we fixe $s,t>0$ and assume that $s=\min(s,t)>0$.
For the lower bound, we have for any fixed $0<a<1$,
since $\psi$ is non-increasing,
$$
 \psi(s/r)+\psi(t/r)\leq 2 \psi(s/r)\leq \psi(as/r)
 $$
 for any $r>0$ such that $s/r\geq s_a$.
Again by monotonicity of   $\psi$,
 $$
as/r\leq  \psi^{-1}\left( \psi(s/r)+\psi(t/r)\right).
$$
So, for any   $0<a<1$, there exists $r_a>0$ such that, for any $0<r<r_a$,
$$
as\leq G_r(s,t).
$$
We deduce, for any $0<a<1$,
$$
a\min(s,t)=as\leq \liminf_{r\rightarrow 0}G_r(s,t)
\leq \limsup_{r\rightarrow 0}G_r(s,t)
\leq \min(s,t)
$$
We get the result with $a\rightarrow 1$. This concludes the proof.
\\

Note that $\theta^{(1)}$ doesn't satisfy the condition $(H_a)$ for any $0<a<1$. We have seen above that $\lim_{r\rightarrow0}G_r(s,t)=\frac{st}{s+t},\; (s,t>0)$
which is strictly less than $\min(s,t)$.
\\

We denote by $f(r):=G_r(s,t)$ with fixed $s,t\in\R$.
We want to prove that $G_0(s,t):=\lim_{r\rightarrow 0^+}f(r)$ exists under some natural condition on $\psi$ with $s,t>0$. This existence is insured if $f(r)$ is non-increasing on some interval
$(0,\varepsilon)$ that is $f'(r)\leq 0$. We give a necessary and sufficient  condition on $\psi$ to fulfil this last condition.
Let $\psi$ as above $(\psi=1-\theta$) and let $V:=(-\psi'\,o\, \psi^{-1})\times\psi^{-1}$. We say that $V$ is locally sub-additive at $0^+$ if there exists $\eta>0$ such that, for all $0<\alpha,\,\beta,\, 
\alpha+\beta<\eta$,  we have
$$
V(\alpha+\beta)\leq V(\alpha)+V(\beta).
$$

We can express the fact that $f'(r)\leq 0 $ is equivalent to this property on $V$.

\begin{Th}\label{limexist}
Let $\psi=1-\theta$ with $\theta :\R\rightarrow (-\infty,1)$ of class $C^2$ such that $\theta'>0$ and 
$\theta''\leq 0$ (i.e. $\psi$ convex). Suppose that $V:=(-\psi'\,o\, \psi^{-1})\times \psi^{-1}$ is locally sub-addditive at $0^+$. Then for any $s,t>0$, 
 $G_0(s,t):=\lim_{r\rightarrow 0^+}f(r)$ exists and 
 $G_0(s,t)\leq \min(s,t)$.
 Moreover, if
there exists $r_0>0$ such that $f'\leq 0$ 
and
$rf'(r)\leq f(r)- f(0^+), 0<r\leq r_0$ then, for any $r\in(0,r_0)$,
\begin{equation}\label{speed}
-r\,\frac{\left(f(0^+)-f(r_0)\right)}{r_0}+f(0^+)\leq f(r)\leq f(0^+).
\end{equation}
\end{Th}

{\bf Comment}:
The   bounds of $f(r)$ in (\ref{speed}) give useful information for numerical simulation because of the bounds 
$0\leq   f(0^+)-f(r) \leq r\,\frac{\left(f(0^+)-f(r_0)\right)}{r_0}
\leq
r\,\frac{\left(\min(s,t)-f(r_0)\right)}{r_0}
$.
\\

\proof 

 Let $f(r):=G_r(s,t)$ with $s,t>0$ and $r>0$.
 Let $H:=-\psi'=\theta'>0$ and $H'=\theta''\leq 0$. So $H$ is positive and non-increasing. A simple computation gives us
 $$
 rf'(r)= f(r)-\frac{sH(s/r)+tH(t/r)}{(H\,o\,\psi^{-1})(\psi(s/r)+\psi(t/r))}.
 $$
 The condition $f'(r)\leq 0$  is equivalent 
 to
 $$
  \psi^{-1}(\psi(s/r)+\psi(t/r))
\leq\frac{\frac{s}{r}H(\frac{s}{r})+\frac{t}{r}H(\frac{t}{r})}{(H\,o\,\psi^{-1})(\psi(s/r)+\psi(t/r))}.
 $$
 Let 
 $\alpha=\psi(s/r)$ and $\beta=\psi(t/r)$.
 Then $f'(r)\leq 0$ if and only if
 $$
  \psi^{-1}(\alpha+\beta)
\leq
\frac{\psi^{-1}(\alpha)\,(H\,o\, \psi^{-1})(\alpha)+\psi^{-1}(\beta)\, (H\,o\, \psi^{-1})(\beta))}{(H\,o\,\psi^{-1})( \alpha+\beta)}.
 $$ 
 Because $H>0$,
 this condition is exactly the sub-additivity property
 \begin{equation}\label {vee}
 V(\alpha+\beta)\leq V(\alpha)+V(\beta).
 \end{equation}
 
 Now assume that $V$  is locally sub-addditive at $0^+$ that is there exists $\eta>0$ such that   
 for all $0<\alpha,\beta,\alpha+\beta<\eta$,  we have
$
V(\alpha+\beta)\leq V(\alpha)+V(\beta).
$
 Fix $s,t>0$ and let $\varepsilon>0$ such that 
 $0<  \max(\psi(s/\varepsilon),\psi(t/\varepsilon))
 \leq \eta$.
 This is possible because $\psi(+\infty)=0^+$
 and 
 $\psi>0$ ($\theta<1$).
  For $r\in (0,\varepsilon)$, we have 
 $0<\alpha:=\psi(s/r)<\eta $ and $0<\beta:=\psi(t/r)<\eta$ since $\psi$ is decreasing. Hence,
 $V(\alpha+\beta)\leq V(\alpha)+V(\beta)$ which implies $f'(r)\leq 0$ for any $0<r< \varepsilon$.
 As a consequence $f(0^+):=\lim_{r\rightarrow 0^+} f(r)$ exists because
 $f(r)$ is always bounded by $\min(s,t)$ (see Proof of  Thm. \ref{zero}).
 \\
 
 Now, we prove (\ref{speed}).
 We have assumed that 
 $$
 rf'(r)\leq f(r)-f(0^+),\; 0<r \leq r_0.
 $$
 So,
 $$
 \left(\frac{f(r)}{r}\right)'= \frac{rf'(r)-f(r)}{r^2}\leq \frac{-1}{r^2}f(0^+).
 $$
Let $0<t<r_0$, by integration on $[t,r_0]$ of the inequality just above, we get
 $$
 \frac{f(r_0)}{r_0} -\frac{f(t)}{t}\leq f(0^+)\left( \frac{1}{r_0}-\frac{1}{t}\right).
 $$
 This can be written as
 $$
 -t\left[\frac{f(0^+)-f(r_0)}{r_0}\right]+f(0)\leq f(t),\quad  0<t\leq r_0.
 $$
 This is the desired result and  completes the proof.
\\

{\bf Two examples}:
In both examples below, it is easy to check that $V$ is sub-addditive on $(0,+\infty)$.
\begin{enumerate}
\item
$\psi^{(1)}(x)=\frac{1}{x+1}, x\geq 1$. So,
$V(y)=y-y^2, \, 0<y<1$ or $V(y)=1-y, \,y>1$.  
\item
$\psi^{(2)}(x)=e^{-x},x\in \R$. So,  $V(y)=-y\ln y$, $0<y<\infty$. \end{enumerate}

The condition  $rf'(r)\leq f(r)- f(0) $  with $0<r$ small enough is satisfied by these two examples.
Let $K=H\,o\,\psi^{-1}=-\psi'\,o\,\psi^{-1}$. With the notations above, we have for $\alpha,\beta>0$:
$$
rf'(r)=f(r)-\frac{sK(\alpha)+tK(\beta)}{K(\alpha+\beta)}.
$$
\begin{enumerate}
\item
For $\psi^{(1)}$: 
$$
\frac{sK(\alpha)+tK(\beta)}{K(\alpha+\beta)}=s\left(\frac{\alpha }{\alpha+\beta}\right)^2
+ t\left(\frac{\beta}{\alpha+\beta}\right)^2
\geq
\inf_{0\leq \lambda\leq 1} \{s\lambda^2+t(1-\lambda)^2\}
=\frac{st}{s+t}=f(0).
$$
In that case, we have equality:
$$
rf'(r)=f(r)-f(0).
$$

\item
For $\psi^{(2)}$: 
Since $K=Id$,
$$
\frac{sK(\alpha)+tK(\beta)}{K(\alpha+\beta)}=\frac{s\alpha+t\beta}{\alpha+\beta}
\geq \min(s,t)=f(0).
$$
Thus  $rf'(r)\leq f(r)- f(0) $.

\end{enumerate}

Now, we give a necessary and sufficient condition on $G(s,t)=\psi^{-1}\left(\psi(s)+\psi(t)\right)$ to be concave
in $(s,t)$ with $s,t>0$.
First, note that $G$  is concave iff $G_r(s,t)=r\,G(s/r,t/r)$ is concave.
($G=G_1$ with this notation).

\begin{Th}\label{concave}
Assume that $\psi: \R \rightarrow (0,+\infty)$ satisfies $\psi'<0$, $\psi''>0$ (i.e. convex).
Let 
$$
L(\alpha):=-\frac{(\psi'\,o\,\psi^{-1})^2}{\psi'' \,o\,\psi^{-1}}(\alpha),\quad \alpha \in \psi(\R).
$$

The following statements are equivalent:
\begin{enumerate}
\item
$G$ is concave in the argument $(s,t)$.
\item
$L$ is non-increasing and sub-additive i.e.
$$
L(\alpha+\beta)\leq L(\alpha)+L(\beta),\quad \alpha,\beta \in \psi(\R).
$$
\end{enumerate}
\end{Th}

Note that $\psi(\R)$ the image of $\R$ by $\psi$ is a subset of $(0,+\infty)$.
\\

\proof
To simplify the presentation of our results, we denote
by $\alpha=\psi(s)$ and $\beta=\psi(t)$,
$$
W:=W(\alpha+\beta)=(\psi'\,o\,\psi^{-1})(\alpha+\beta)<0, \quad \alpha+\beta \in \psi(\R),
$$ 
and
$$
U:=U(\alpha+\beta)=(\psi''\,o\,\psi^{-1})(\alpha+\beta)>0, \quad \alpha+\beta \in \psi(\R).
$$

A rather boring computation gives us, 
$$
R:={\partial}_{s,s}G(s,t)=
\left[ \psi''(s) W^2-(\psi'(s))^2U\right]/W^3,
$$

$$
T:= {\partial}_{t,t}G(s,t)=
\left[ \psi''(t) W^2-(\psi'(t))^2U\right]/W^3,
$$

$$
S:= {\partial}_{s,t}G(s,t)=
-\psi'(s)\psi'(t)\frac{U}{W^3}.
$$
It is well-known that $G$ is concave iff $R\leq 0, T\leq 0$ and $RT-S^2\geq 0$.
For  the condition $R\leq 0$ (similarily for $T\leq 0$),   we get due to the fact that $W^3< 0$,
 $$
 \psi''(s) W^2\geq
 (\psi'(s))^2U.
$$
This can written as
$$
- \frac{(\psi'(s))^2}{ \psi''(s)}\geq  -\frac{W^2}{U}.
 $$
That is, with $s=\psi^{-1}(\alpha)$ and $t=\psi^{-1}(\beta)$,
$$
L(\alpha)\geq L(\alpha+\beta)
$$
 This expresses  the fact that $L$ is non-increasing.
\\

For  the condition $RT-S^2\geq 0$,   we obtain after simplifications:
$$
W^6(RT-S^2)=
\psi''(s)\,\psi''(t)W^4-
\left[
\psi''(s)(\psi')^2(t)+\psi''(t)(\psi')^2(s)\right]W^2U.
$$
By  similar manipulations as in the case   $R\leq 0$, we can express the condition
$RT-S^2\geq 0$ by
$$
L(\alpha+\beta)\leq L(\alpha)+L(\beta),\quad \alpha,\beta,\alpha+\beta \in \psi(\R).
$$
that is $L$ is sub-additive. The proof is completed.
\\

Note that $T\leq 0$ iff $R\leq 0$. Indeed,  $R\leq 0$ iff $L$ is non-increasing. Thus, by symmetry arguments in $(s,t)$, 
 $T\leq 0$ iff $L$ is non-increasing. 
 \\

An examination of the two examples above leads to the  following  important remark. We compute the functions $L$ for the examples  
  $\psi^{(1)}(x)=\frac{1}{x+1}$ and $\psi^{(2)}(x)=e^{-x}$ and  obtain  $L_1(\alpha)=-\frac{1}{2} \alpha$ and $L_2(\alpha)=-\alpha$. They are additive functions !
  It suggests to find a one parameter family of function $\phi_{\lambda}$ 
  giving as particular cases the functions  
  $\psi^{(1)}$ and $\psi^{(2)}$.
  This can be done by solving the equation $RT-S^2=0$. Indeed, additivity of $L$ exactly correspond to this equation. The equation $RT-S^2=0$ can be written as
  $L(\alpha)=-\frac{1}{\lambda}\alpha,\, \lambda>0$ ($L$ is non-increasing) or equivalently 
  $$
  (\psi')^2=\frac{1}{\lambda}\psi\,\psi''.
  $$
First case: $\lambda>1$. We obtain as solution of this equation
$$
\phi_{\lambda}(x)=\frac{1}{(c_1x+1)^{\frac{1}{\lambda -1}}}, \, \quad x>\frac{-1}{c_1},
$$
for some $ c_1 >0$.
\\

Second case $\lambda=1$: this case  has to be solved independently,
we get
$$
\phi_{1}(x)=e^{-Dx},\, x\in \R,
$$
for some $D>0$.
This case can be seen as a limit case as $\lambda \rightarrow 1^+$. The function $\phi_{1}$ is different of nature
of the functions $\phi_{\lambda}, \lambda>1$.
\\

The condition $RT-S^2=0$ introduces a family of new examples of $\psi$ satisfying the conditions
of Theorem  \ref{concave} namely 
$\phi_{\lambda}$. This condition is exactly
$\det \,(Hess G)=0$. That is the Hessian has an eigenvalue $\mu_1=0$ and $\mu_2\leq0$ (since the trace
of the Hessian matrix is $R+T\leq 0$).

\section{Convergence and error estimate}

Let $H_r(x):=G_r(x,F(x))=\left( G_r(x_i,F_i(x))\right)_{i=1}^n $ with $G_r$ defined as above.
When $F$ is a $(P_0)$-function, we can easily prove the following result.
\begin{lem}
Assume that $F$ is a $(P_0)$-function then $H_r$ is  
$(P)$-function for any $r>0$.
\end{lem}
\proof For any $x,y\in \R^n$, there exits $i:1...n$ such that $x_i\neq y_i$. We can assume that 
$x_i> y_i$ and $F_i(x)\geq F_i(y)$. Since $\psi$ is decreasing
:
$\psi(x_i/r)< \psi(y_i/r)$ 
and 
$ \psi(F_i(x)/r)\leq \psi(F_i(y)/r)$.
Consequently,
$\psi(x_i/r)+ \psi(F_i(x)/r)< \psi(y_i/r)+\psi(F_i(y)/r)$.
Again by the monotonicity of $\psi^{-1}$, 
$G_r(x_i,F_i(x))> G_r(y_i,F_i(y))$.
 Hence, $H_r$ is  
$(P)$-function for any $r>0$.

\noindent So that, when  the assumptions of Theorem \ref{limexist} are satisfied, we obtain the following convergence result.

\begin{Th}
Under the assumptions of Theorem \ref{limexist}, assume that $F$ is a $(P_0)$-function and that  $(H1)$ is  satisfied (The solution set of the NCP is nonempty and compact).\\
(i) There exists an $\hat r>0$ such that for any $0<r<\hat r$,  $H_r(x)=0$ will have a unique solution $x^{(r)}$, the mapping $r\to x^{(r)}$ is continuous on $(0,\hat r)$,  and \\
(ii) $ \displaystyle \lim_{r\to 0} dist (x^{(r)},$ \cal{Z} $) =0$.
\end{Th}
\proof
 This is a dierct application of (\cite{Se-Ta} Theorem 4 (2)).\\
\begin{rem}
Under  the assumptions of Theorem \ref{concave}, we can prove an other convergence result based on the smoothing
thechnique disscussed in {\rm \cite{BT}}.
\end{rem}

\noindent When using $\theta \ge \theta^{(1)}$ on $\R^+$ (this is the case of $\theta^{(2)} $for example), we can prove  an estimate for  the error term $ \vert\vert x^*-x^{(r)}\vert\vert $ between the solution $x^*$ and the approximation $x^{(r)}$ under an assumption of monotonicity of $F$.

\begin{prop}
Assume that  $\theta \ge \theta^{(1)}$ on $\R^+$, that $x^*$ is a solution of $<x^*,F(x^*)>=0$ and 
$x^{(r)}$ is a nonnegative solution of $H_r(x)=0$.\\
(i) $x_i^{(r)}F_i(x^{(r)}) \le r^2$  $\forall i=1\dots n$\\
(ii) If $F$ satisfies the  condition: 
$$
h(\vert\vert x-y\vert\vert) \leq\; <x-y,F(x)-F(y)> 
$$
with  $h:\R^+\longrightarrow \R$ such that 
$h(0)=0$
and there exists
  $\varepsilon>0$ such that 
$h:[0,\varepsilon)\longrightarrow [0,\eta)$  
is  an increasing bijection.
Then there exists $r_0>0$ such that for any $r\in (0,r_0)$,
\begin{equation}\label{error}
 \vert\vert x^*-x^{(r)}\vert\vert
 \leq h^{-1}(n r^2).
\end{equation}
\end{prop}

The inequality  (\ref{error}) gives the maximal behavior of the error
 in terms of the function $h$.
\\

\proof
(i) $x^{(r)}$ satisfies  $H_r(x^{(r)})=0$, so that 
$$
\theta(\frac{x_i^{(r)}}{r})+\theta(\frac{F_i(x^{(r)}}{r})=1, \qquad i:1...n.
$$
and since $\theta \ge \theta^{(1)}$, we obtain 
$$
\theta^{(1)}(\frac{x_i^{(r)}}{r})+\theta^{(1)}(\frac{F_i(x^{(r)}}{r})\leq 1, \qquad i:1...n.
$$
Then, a simple calculus yields to 
$$x_i^{(r)}F_i(x^{(r)}) \le r^2 \qquad i:1...n.$$

(ii) We have 
$$
<x^*-x^{(r)},F(x^*)-F(x^{(r)})>=
$$
$$
<x^*,F(x^*)> - <x^*,F(x^{(r)})>-
<(x^{(r)},F(x^*)>+<x^{(r)},F(x^{(r)})>\; 
\leq n r^2.
$$
Indeed,  the first term of the R-H-S is zero, the two middle terms are non-positive and the last term is $n r^2$.
So, by monotonicity,
$$
 h(\vert\vert x^*-x^{(r)}\vert\vert) 
\leq  n r^2
 $$
  Let $r_0$ such that $r^2_0 <\eta$.
Since 
 $h$ is a bijection from $[0,\varepsilon)$ onto $[0,\eta)$
 and $h^{-1}$ is increasing, we obtain
 $$
  \vert\vert x^*-x^{(r)}\vert\vert
 \leq h^{-1}(n r^2).
 $$
 The proof is completed.


\section{Numerical results}
In order to verify the theoretical assertions, we present some numerical experiments for  two smoothing approaches using the $\theta$ functions $\theta^{(1)}$ and $\theta^{(2)}$.  We first consider a simple 2-dimensional problem which is analytically solvable 
where
$$F(x,y)=(2-x-x^3,y+y^3-2)^T.$$
The unique optimal solution for this problem is $(0,1)$. \\
The following figure presents the evolution of the second coordinate of the iterates for the two smoothing functions. The optimal solution was reached up to a tolearance of $10^{-10}$ in no more than 3 iterations. We use the updating strategy for the penalization parameter as precised later. The red and green points correspond respectively to the iterates of the smoothing method
$\theta^{(1)}$ and $\theta^{(2)}$.\\
 \begin{figure}[ht!]
\begin{center}
{\label{gout}\includegraphics[width=8cm]{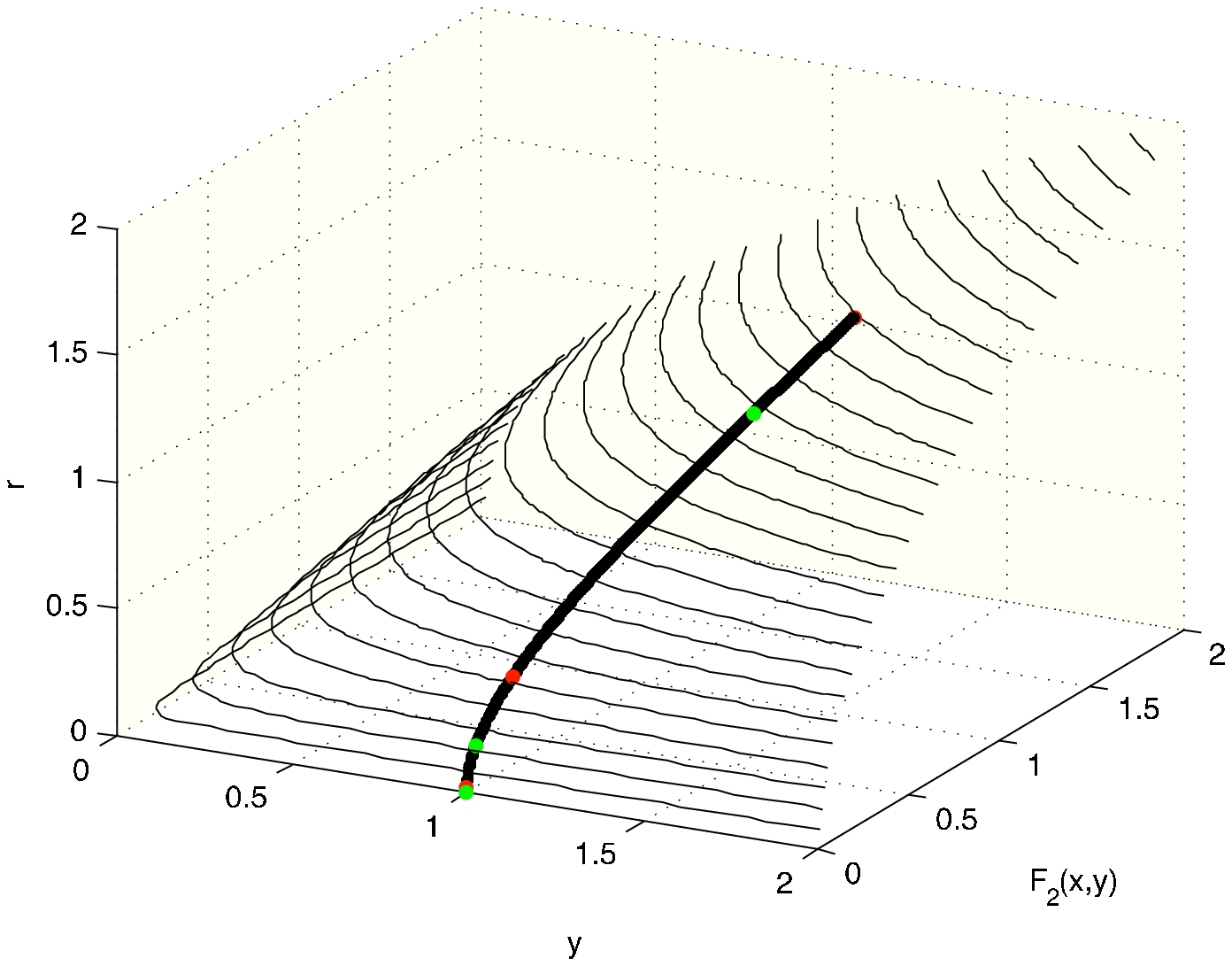}}
\caption{Evolution of $(y,F_2(x,y),r)$.}\label{figcomp}
\end{center}
\end{figure}

 
We also consider a set of 10 NCP test problems with different and varying number of variables. \\
For each test problem and each smoothing function, we use 11 different starting points:  a vector of ones and 10 uniformly generated vectors with entries in $(0,20)$.\\
The starting value for the smoothing parameter is fixed with respect to the theoretical properties as
$$r^0=\max ( 1, \sqrt{ \max_{1\le i\le n} \vert x^0_i \bullet F_i(x^0)\vert  } ).$$
This parameter is then updated as follows
$$r^{k+1}=\min (0.1r^k,({r^k})^2, \sqrt{\max_{1\le i\le n}\vert x^k_i \bullet F_i(x^k)\vert } )$$
until the stopping rule
$$\max_{1\le i\le n}\vert x^k_i \bullet F_i(x^k)\vert \le 10^{-8}$$
is satisfied.\\

Precise descriptions of the test problems   P1 and P2 can be found in \cite{HW}.
 P3 is a test problem from  \cite{LZ}  while P4 and P5 can be found in  \cite{DY} , these two problems correspond to a non-degenerate and a degenerate examples  of Kojima-Shindo NCP test problems \cite{KS}.  
 The other test problem are described in\cite{Tin,Har}. They correspond respectively to the NASH-COURNOT test problems with $n=5$ and $n=10$ and to the HpHard test problem with $n=20$ , $n=30$ and $n=100$.\\
We used a standard laptop (2.5 Ghz, 2Go M) and a very simple matlab program using the $fsolve$ function. \\

We list  in the following table, the worst obtained results.  $n$ stands for the number of variables. OutIter is the number of outer iterations (number of changes of the smoothing parameter $r$) and InIter corresponds to the total number of jacobian evaluations. Res. and Feas. correspond to the following optimality and feasibility measures
$$\textrm{Res.}=\max_{1\le i\le n}\vert x_i \bullet F_i(x)\vert $$
and 
$$\textrm{Feas.}=\Vert\min(x,0)\Vert_1+\Vert\min(F(x),0)\Vert_1.$$
The results show that the second smoothing function is much more efficient and powerful. This was foreseeable since 
$$\forall x\ge 0 \qquad 1-\delta_0 (x)\ge \theta^{(2)} (x)\ge \theta^{(1)}(x).$$
\newpage

{
\begin{table}[htb]
\begin{displaymath}
\begin{array}{llclclclclclclcll}
\hline\hline
\textrm{Pb} & \textrm{size}&\textrm{OutIter} & \textrm{InIter} & \textrm{Res.} &  \textrm{Feas.} &\textrm{cpu time (s)}\\ 
 &&(\theta_1,\theta_2) &(\theta_1,\theta_2) & (\theta_1,\theta_2) & (\theta_1,\theta_2) &(\theta_1,\theta_2) \\ \hline
\hline \textrm{P1}&\textrm{10}&(6 ,4)&(65 ,15)&(5.6e\!\!-\!\!15 ,2.5e\!\!-\!\!18)&(  1.1e\!\!-\!\!11  ,1.3e\!\!-\!\!10)&(0.22,0.09)\\ 
\hline \textrm{}&\textrm{100}&(6 ,4)&(68 ,19)&(1.6e\!\!-\!\!14 ,7.1e\!\!-\!\!22)&(  5.1e\!\!-\!\!13  ,1.4e\!\!-\!\!14)&(3.73,1.19)\\ 
\hline \textrm{}&\textrm{500}&(6 ,4)&(83 ,21)&(5.4e\!\!-\!\!12 ,1.6e\!\!-\!\!16)&(  1.9e\!\!-\!\!16  ,1.4e\!\!-\!\!14)&(31.15,89.26,)\\ 
\hline \textrm{}&\textrm{1000}&(6 ,5)&(77 ,40)&(3.0e\!\!-\!\!14 ,3.1e\!\!-\!\!14)&(  5.1e\!\!-\!\!18  ,1.8e\!\!-\!\!17)&(388.59,201.43)\\ \hline
\hline \textrm{P2}&\textrm{10}&(6,4)&(79,23)&(2.1e\!\!-\!\!15,2.7e\!\!-\!\!15)& (7.6e\!\!-\!\!11, 9.6e\!\!-\!\!19)  &(0.31,0.11)\\ 
\hline \textrm{}&\textrm{100}&(6,4)&(88,33)&(1.84e\!\!-\!\!12,1.0e\!\!-\!\!23)&   (7.1e\!\!-\!\!10,3.1e\!\!-\!\!14)&(4.83,1.80)\\ 
\hline \textrm{}&\textrm{500}&(6,4)&(96,41)&(6.5e\!\!-\!\!10,1.9e\!\!-\!\!16)&  (6.6e\!\!-\!\!09,1.2e\!\!-\!\!12)  &(112.14,49.59)\\ 
\hline \textrm{}&\textrm{1000}&(6,5)&(114,67)&(1.0e\!\!-\!\!17,1.4e\!\!-\!\!23)&( 2.4e\!\!-\!\!08 ,7.5e\!\!-\!\!18)&(530.42,328.15)\\ \hline
\hline \textrm{P3}&\textrm{10}&(5,4)&(63,15  ) & ( 2.2e\!\!-\!\!12,2.7e\!\!-\!\!21)&(4.9e\!\!-\!\!08,1.4e\!\!-\!\!11 )   &(0.22,0.09)\\ 
\hline \textrm{}&\textrm{100}&(5,4)&(71,18  )  &(7.9e\!\!-\!\!13, 2.6e\!\!-\!\!15)& (9.5e\!\!-\!\!08,4.5e\!\!-\!\!08 ) &(3.10,1.02)\\ 
\hline \textrm{}&\textrm{500}&(5,4)&(73,21 )   &(1.1e\!\!-\!\!14,2.6e\!\!-\!\!16)& (1.5e\!\!-\!\!07, 5.9e\!\!-\!\!09 ) &(78.11,26.15)\\ 
\hline \textrm{}&\textrm{1000}&(5,4)&(81,26)&(6.1e\!\!-\!\!13,1.2e\!\!-\!\!15) &( 8.2e\!\!-\!\!10      ,2.4e\!\!-\!\!16)& (335.37,138.23)\\ \hline
\hline \textrm{P4}&\textrm{4}&(6,4)&(63,20)&(5.4e\!\!-\!\!12,3.2e\!\!-\!\!17)&(6.1e\!\!-\!\!09,2.8e\!\!-\!\!12)&(0.15,0.08)\\ \hline
\hline \textrm{P5}&\textrm{4}& (6,4)&(141,23)& (9.8e\!\!-\!\!14,2.1e\!\!-\!\!23)&(3.4e\!\!-\!\!07,3.2e\!\!-\!\!12)&(0.28,0.06)\\  \hline
\hline \textrm{P6}&\textrm{5}&(5,3)&(47,17)& (1.3e\!\!-\!\!14, 4.3e\!\!-\!\!27)&(4.9e\!\!-\!\!12,8.1e\!\!-\!\!17)&(0.16,0.07)\\ \hline
\hline \textrm{P7}&\textrm{10}&(6,4)&(110,33)& (1.2e\!\!-\!\!16, 6.1e\!\!-\!\!19)&(1.1e\!\!-\!\!12,4.5e\!\!-\!\!14)&(0.37,0.14)\\  \hline
\hline \textrm{P8}&\textrm{20}&(6,5)&(145,66)& (2.9e\!\!-\!\!13, 3.7e\!\!-\!\!21)&(0,4.4e\!\!-\!\!12)&(1.33,0.46)\\  
\hline \textrm{P9}&\textrm{30}&(6,6)&(106,77)& (3.7e\!\!-\!\!14, 9.6e\!\!-\!\!21)&(4.4e\!\!-\!\!08,6.4e\!\!-\!\!11)&(2.24,0.85)\\  
\hline \textrm{P10}&\textrm{100}&(6,6)&(209,113)& (8.5e\!\!-\!\!11, 2.1e\!\!-\!\!23)&(2.1e\!\!-\!\!07,1.8e\!\!-\!\!12)&(42.09,19.12)\\  \hline \hline

\end{array}
\end{displaymath} \caption{Results for  $\theta^{1}$ and  $\theta^{2}$}\label{kk}
\end{table}
}



\begin{thebibliography}{BCLS}

\bibitem[ACH]{ACH} 
 Auslender, A.; Cominetti, R.; Haddou, M. Asymptotic analysis for penalty and barrier methods in convex and linear programming. Math. Oper. Res. 22 (1997), no. 1, 43--62. 

\bibitem[BT]{BT} 
Ben-Tal, A. and M.Teboulle. A Smoothing Technique for Nondifferentiable Optimization 
Problems. In Dolecki, editor, Optimization, Lectures notes in Mathematics 1405, pages 1–11, 
New York, 1989. Springer Verlag. 

\bibitem[DY]{DY} 
Ding, Jundi; Yin, Hongyou A new homotopy method for nonlinear complementarity problems.  Numer. Math. J. Chin. Univ. (Engl. Ser.)  16  (2007),  no. 2, 155--163.

\bibitem[FMP]{FMP} 
Complementarity: applications, algorithms and extensions. Papers from the International Conference on Complementarity (ICCP99) held in Madison, WI, June 9--12, 1999. Edited by Michael C. Ferris, Olvi L. Mangasarian and Jong-Shi Pang. Applied Optimization, 50. Kluwer Academic Publishers, Dordrecht, 2001.

\bibitem[FP]{FP} 
Ferris, M. C.; Pang, J. S. Engineering and economic applications of complementarity problems. 
SIAM Rev. 39 (1997), no. 4, 669--713.

 \bibitem[Had]{Had} 
 Haddou M. A new class of smoothing methods for mathematical programs with equilibrium 
constraints. 
Paciﬁc Journal of Optimization, vol 5(1) (2009) , pp.86-96. 

 \bibitem[Har]{Har} 
P.T. Harker. Accelerating the convergence of the diagonalization 
 and projection algorithms for finite-dimensional variational inequalities'
Mathematical Programming 48, (1990) pp. 29-59.

 \bibitem[HW]{HW} 
 HUANG C,  WANG S. A power penalty approach to a Nonlinear Complementarity Problem.
Operations research letters 2010, vol. 38, no1, pp. 72-76.

\bibitem[KS]{KS} 
 Kojima M, Shindo S. Extensions of Newton and quasi-Newton methods to systems of PC1 
equations. J. Oper. Res. Soc. Jpn., 1986, 29: 352-374. 

 
\bibitem[LZ]{LZ} 
 Dong-hui Li, Jin-ping Zeng.
 A penalty technique for nonlineair problems.
  Journal of Computational Mathematics, Vol.16, No.1, 1998, 40–50.

\bibitem[PL]{PL} 
Peng, J.M. and Z. Lin. A Non-interior Continuation Method for Generalized Linear Complementarity Problems. Mathematical Programming, 86:533–563, 1999.

\bibitem[Se-Ta]{Se-Ta} 
Seetharama Gowda, M.; Tawhid, M. A. Existence and limiting behavior of trajectories associated with $P_0$-equations. Computational optimization---a tribute to Olvi Mangasarian, Part I.  Comput. Optim. Appl.  12  (1999),  no. 1-3, 229--251.



\bibitem[Tin]{Tin} 
{\verb" http://dm.unife.it/pn2o/software/Extragradient/test_problems.html"}

\end{thebibliography}
\end{document}